\documentclass[11pt]{amsart}

\usepackage[letterpaper,centering,margin=1.47in]{geometry}

\usepackage{amsthm}
\usepackage{amssymb}
\usepackage{amsmath}
\usepackage{amsxtra}
\usepackage[all,cmtip]{xy}

\title{The Atiyah--Segal completion theorem in twisted $K$-theory}

\author{Anssi Lahtinen}
\email{lahtinen@math.ku.dk}
\urladdr{http://www.math.ku.dk/~lahtinen/}

\newtheorem{theorem}{Theorem}
\newtheorem{proposition}{Proposition}
\newtheorem{lemma}{Lemma}

\newtheorem{corollary}{Corollary}
\theoremstyle{remark}
\newtheorem{remark}{Remark}

\makeatletter 
\let\c@proposition=\c@theorem
\let\c@lemma=\c@theorem
\let\c@conjecture=\c@theorem
\let\c@corollary=\c@theorem
\let\c@remark=\c@theorem
\makeatother

\newcommand{\Z}{\mathbb{Z}}
\newcommand{\R}{\mathbb{R}}

\newcommand{\T}{\mathbb{T}} 
\newcommand{\pt}{\mathrm{pt}}
\newcommand{\isom}{\approx}
\newcommand{\smashprod}{\wedge}
\newcommand{\invlim}{\varprojlim}
\newcommand{\incl}{\hookrightarrow}

\newcommand{\compl}{\sphat}
\newcommand{\Vect}{\mathrm{Vect}}
\newcommand{\reduced}[1]{\tilde{#1}}
\newcommand{\rK}{\reduced{K}}
\newcommand{\proK}{\mathbf{K}}
\newcommand{\rproK}{\reduced{\proK}}
\newcommand{\tH}{{\tilde H}}
\newcommand{\tG}{{\tilde G}}
\newcommand{\quotb}[2]{#1 / #2 \cdot #1}
\newcommand{\modmod}{{\slash\!\!\slash}}
\newcommand{\quotient}{\, / \,}
\DeclareMathOperator{\colim}{colim}

\begin{document}

\begin{abstract}   

A basic result in equivariant $K$-theory, the Atiyah--Segal 
completion theorem relates the $G$-equivariant $K$-theory of a 
finite $G$-CW complex to the non-equivariant 
$K$-theory of its Borel construction.
We prove the analogous result for twisted equivariant $K$-theory.
\end{abstract}

\maketitle

\section{Introduction}

The aim of this note is to prove the following twisted analogue of the
Atiyah--Segal completion theorem \cite{ASC}.
\begin{theorem}
\label{mainthm}
Let $X$ be a finite $G$-CW complex, where $G$ is a compact Lie group.
Then the projection map $\pi : EG \times X \to X$ induces an isomorphism
\begin{displaymath}
  K^{\tau+*}_G(X)\compl_{I_G} \xrightarrow{\,\isom\,}
  K^{\pi^*(\tau)+*}_G(EG \times X)
\end{displaymath}
for any twisting $\tau$ corresponding to an element of
$H^1_G(X;\,\Z/2) \oplus H^3_G(X;\, \Z)$.
\end{theorem}
\noindent
Here $I_G \subset R(G)$ is the augmentation ideal of the
representation ring $R(G)$ and $(-)\sphat_{I_G}$ indicates completion.
The classical theorem is the case $\tau = 0$. 
As in the untwisted case, Theorem \ref{mainthm} 
implies a comparison between the equivariant $K$-theory
of a finite $G$-CW complex $X$ and the 
non-equivariant $K$-theory of its Borel construction.
\begin{corollary}
Let $X$ be a finite $G$-CW complex. Then there exists an isomorphism
\[
  K^{\tau+*}_G(X)\compl_{I_G} 
  \isom   K^{\tau+*}(EG \times_G X)
\]
for any twisting $\tau$ corresponding to an element of 
\[
H^1_G(X;\,\Z/2) \oplus H^3_G(X;\, \Z) 
= 
H^1(EG \times_G X;\,\Z/2) \oplus H^3(EG \times_G X;\, \Z). \qed
\]
\end{corollary}
Theorem \ref{mainthm}
generalizes a result by C~Dwyer, who has proved the theorem in the
case where $G$ is finite and the twisting $\tau$ is given by a cocycle
of $G$~\cite{CD}. While versions of the theorem for compact Lie groups
have been known to experts (for example, such a theorem is used in the
proof of \cite{FHTcompl}, Proposition 3.11), to our knowledge no proof
of the general theorem appears in the current literature. Our goal is
to fill in this gap.

We shall prove Theorem \ref{mainthm} in two stages. 
First we prove the theorem in the special case of a twisting arising from a 
graded central extension 
\[
    1\to \T \to G^{\tau} \to G \to 1,\qquad \epsilon : G \to \Z/2
\]
of $G$ by the circle group $\T$.
For such twistings, twisted $G$-equivariant $K$-groups correspond to
certain direct summands of untwisted ${G^\tau}$-equivariant $K$-groups, and
the Adams--Haeberly--Jackowski--May argument contained in \cite{AHJM1}
goes through with these summands to prove the theorem in this case.
It follows that the theorem holds when $X$
is a point, and the general theorem then follows by a Mayer--Vietoris
argument.

As our definition of twisted $K$-theory, we use
Freed, Hopkins and Teleman's elaboration \cite{FHT1} of the
Atiyah--Segal model developed in \cite{AStwistedK}. Thus for a
$G$-space $X$, the notation $K^{\tau + *}_G(X)$ is a shorthand for
$K^{\tau+*}(X\modmod G)$, where $X\modmod G$ is the quotient groupoid of $X$.
Of course, the completion theorem should remain true in any 
reasonable alternative model for twisted equivariant $K$-theory as well.

This note is structured as follows. In section~\ref{coh} we describe 
a pro-group valued variant of $K$-theory which we shall employ in 
section~\ref{sc} to handle the case of a twisting given by a central
extension. Section~\ref{gc} then contains a proof of the general theorem.

\subsection*{Acknowledgements}

The present research formed a part of my Stanford University
PhD thesis, and I would like to express my warm thanks to my 
PhD advisor, Professor Ralph Cohen, for his advice and support. 
In addition, it is my pleasure to acknowledge support from 
the Finnish Cultural Foundation and from 
the Danish National Research Foundation (DNRF) through the 
Centre for Symmetry and Deformation.

\section{A convenient cohomology theory}
\label{coh}

We shall now describe a cohomology theory which will be used in the
next section to prove the completion theorem for twistings arising
from a graded central extension. Let 
\[
   1 \to C \to \tG \to G \to 1 
\]
be a central extension of a compact Lie group $G$ by a commutative 
compact Lie
group $C$, and let $X$ be a finite $G$-CW complex. Via the map $\tG \to G$
we can  view $X$ as a $\tG$-space on which $C$ acts trivially. The
semigroup $\Vect_\tG(X)$ of isomorphism classes of $\tG$-equivariant
vector bundles over $X$ decomposes as a direct sum
\begin{equation}
\label{vdecomp}
  \Vect_\tG(X) = \bigoplus_{\pi \in \hat{C}} \Vect_\tG(X)(\pi)
\end{equation}
where $\hat{C}$ denotes the set of isomorphism classes of irreducible
representations of $C$, and where $\Vect_\tG(X)(\pi)$ is the semigroup
of isomorphism classes of those $\tG$-vector bundles over $X$
whose fibers 
are $\pi$-isotypical as representations of $C$, that is, 
isomorphic to sums of copies of $\pi$. The decomposition
\eqref{vdecomp} leads to a decomposition
\begin{equation}  
\label{kdecomp}
 K_\tG^*(X) = \bigoplus_{\pi \in \hat{C}}K_\tG^*(X)(\pi)
\end{equation}
and similarly for reduced $K$-groups.\footnote{In fact, tensor product
  makes $K_\tG^*(X)$ into a $\hat C$-graded algebra where the modules
  $K_\tG^*(X)(\pi)$ are the homogeneous parts. However, we shall not
  need this graded algebra structure.}  
Here $K_\tG^0(X)(\pi)$ is the Grothendieck group of
$\Vect_\tG(X)(\pi)$, and $K_\tG^q(X)(\pi)$ for non-zero $q$ is defined
by using the suspension isomorphism and the Bott periodicity map.  By
inspection and definition, the decomposition \eqref{kdecomp} is
compatible with $G$-equivariant maps of spaces, with the suspension
isomorphism, with the Thom isomorphism for $G$-equivariant vector
bundles, and, as a special case, with the Bott periodicity map. Thus
for each $\pi \in \hat{C}$, we can view $K^*_\tG(-)(\pi)$ as a
$\Z/2$-graded cohomology theory defined on finite $G$-CW complexes and
taking values in graded $R(G)$-modules.

Although the decomposition \eqref{kdecomp} fails for infinite $X$ in
general, it is possible to extend each one of the theories
$K^*_\tG(-)(\pi)$ to infinite $G$-CW complexes by means of suitable
classifying spaces.  However, since having the theories available for
finite complexes suffices for most of our purposes, we will not
elaborate this point. Instead, we point the
reader to the proof of Proposition 3.5 in \cite{FHT1} for a
description of the appropriate classifying space when $\pi$ is the
defining representation of the circle group $\T$, which is the only case
where we will need to apply $K^*_\tG(-)(\pi)$ to an infinite complex 
in the sequel.

Our interest in the groups $K^*_\tG(X)(\pi)$ is explained by the
following proposition. Recall that a graded central extension of a
group $G$ is a central extension of $G$ together with a homomorphism
from $G$ to $\Z/2$.

\begin{proposition}[A reformulation of Proposition 3.5 of \cite{FHT1}]
\label{propceiso}
Let $G$ be a compact Lie group, let $X$ be a $G$-space, and let $\tau$
be the twisting given by a graded central extension
\[
    1\to \T \to G^{\tau} \to G \to 1,\qquad \epsilon : G \to \Z/2
\]
of $G$ by the circle group $\T$. Let $S^1(\epsilon)$ denote the
one-point compactification of the 1-dimensional representation of $G$
given by $(-1)^{\epsilon}$. Then there is a natural isomorphism
\[
K^{\tau + n}_G(X) \isom \reduced{K}^{n+1}_{G^{\tau}}(X_+ \smashprod S^1(\epsilon))(1)
\]
where ``\,$1$\!'' refers to the defining representation of $\T$. \qed
\end{proposition}

The groups $K^*_\tG(X)(\pi)$ are not what we are going to use in the
next section. Instead, we need pro-group valued versions completed at
the augmentation ideal. (For background material on pro-groups, we
refer the reader to \cite{AHJM2}.) Given an arbitrary $G$-CW complex
$X$ and an irreducible representation $\pi$ of $C$, we let
$\proK_\tG^*(X)(\pi)$ denote the pro-$R(G)$-module
\begin{equation*}
  \proK_\tG^*(X)(\pi) = \{ K_\tG^*(X_{\alpha})(\pi) \}_{\alpha}
\end{equation*}
where $X_\alpha$ runs over all finite $G$-CW subcomplexes of $X$ and
the structure maps of the pro-system are those induced by inclusions
between subcomplexes.  The groups of our interest are then
the pro-$R(G)$-modules
\[
  \proK_\tG^*(X)(\pi)\compl_{I_G}=  \{ K_\tG^*(X_{\alpha})(\pi) \quotient 
  I_G^n \cdot K_\tG^*(X_{\alpha})(\pi) \}_{\alpha,n},
\]
where $X_\alpha$ again runs over the finite $G$-subcomplexes of
$X$, $n$ runs over the natural numbers, and the structure maps of
the pro-system are the evident ones. We think of
$\proK_\tG^*(X)(\pi)\compl_{I_G}$ as the completion of
$\proK_\tG^*(X)(\pi)$ with respect to the augmentation ideal $I_G$.
Reduced variants $\rproK_\tG^*(X)(\pi)$ and $\rproK_\tG^*(X)(\pi)\compl_{I_G}$
for a based $G$-CW complex $X$ are defined in a similar way using the 
reduced groups $\rK^*_\tG(X_\alpha)(\pi)$, where $X_\alpha$ now runs through the
finite $G$-CW subcomplexes of $X$ containing the base point. The
crucial feature of the groups $\proK_\tG^*(X)(\pi)\compl_{I_G}$ for us
is that they form a cohomology theory on the category of $G$-CW
complexes (and therefore, by $G$-CW approximation, on the category of
all $G$-spaces). Phrased in terms of the reduced groups, this means
that the following axioms hold.
\begin{enumerate}
\item (Homotopy invariance) If $X$ and $Y$ are based $G$-CW complexes
  and $f,g : X \to Y$ are homotopic through
  based $G$-equivariant maps, then the induced maps 
\[f^*,g^* :  \rproK_\tG^*(Y)(\pi)\compl_{I_G}  \to \rproK_\tG^*(X)(\pi)\compl_{I_G}\]
 are equal.
\item (Exactness) If $X$ is a based $G$-CW complex and $A$ is a subcomplex 
 of $X$ containing the base point, then the sequence \[ 
    \rproK_\tG^*(X/A)(\pi)\compl_{I_G} \to \rproK_\tG^*(X)(\pi)\compl_{I_G}
    \to \rproK_\tG^*(A)(\pi)\compl_{I_G}
  \] 
  is pro-exact.
\item (Suspension) For each $q$, there exists a natural isomorphism \[
  \Sigma : \rproK_\tG^q(X)(\pi)\compl_{I_G} \isom
  \rproK_\tG^{q+1}(\Sigma X)(\pi)\compl_{I_G}\]
\item (Additivity) If $X$ is the wedge sum of a family $\{X_i\}_{i \in I}$ 
  of based $G$-CW complexes, the inclusions $X_i \incl X$
  induce an isomorphism
   \[
      \rproK_\tG^*(X)(\pi)\compl_{I_G}
      \xrightarrow{\,\isom\,} \prod_{i\in I} \rproK_\tG^*(X_i)(\pi)\compl_{I_G}
   \]
\end{enumerate}
The only difficulties in verifying these properties 
arise from the exactness axiom.
\begin{proposition}
\label{propexact}
  The functor $\rproK^*_\tG(-)(\pi)\compl_{I_G}$ satisfies the exactness axiom.
\begin{proof}[Sketch of proof]
 As in \cite{AHJM2}, because the ring $R(G)$ is Noetherian (see
 \cite{SRRL}, Corollary 3.3), the result follows from the Artin--Rees
 lemma once $\rK^*_\tG(Z)(\pi)$ is known to be finitely generated as
 an $R(G)$-module for any finite based $G$-CW complex $Z$. We shall
 prove that $\rK^*_\tG(Z)(\pi)$ is finitely generated by reduction to
 successively simpler  cases.  Filtering $Z$ by skeleta and
 using the wedge and suspension axioms, we see that it is enough to
 consider the case where $Z = G/{H}_+$ for some closed subgroup ${H}$
 of $G$. Let ${\tH}$ denote the inverse image of $H$ in $\tG$. Then
 ${\tH}$ is a central extension of ${H}$ by $C$, and we have
 $\tG$-equivariant isomorphisms
  \[ 
     G/{H} \isom (\tG/C)/ ({\tH}/C) \isom \tG/{\tH}.
  \]
  The $R(G)$-module isomorphisms 
  \[
      K_\tG^*(G/{H}) \isom K_\tG^*(\tG/{\tH}) \isom
      K_{{\tH}}^*(\pt)
  \]
  preserve the direct sum decomposition \eqref{kdecomp}, whence 
  we obtain an isomorphism 
  \[
      K_\tG^*(G/{H})(\pi) \isom K_{{\tH}}^*(\pt)(\pi).
  \]
  Here the latter group can be identified with the summand
  $R({\tH})(\pi)$ of $R({\tH})$ generated by those representations of
  ${\tH}$ which restrict to $\pi$-isotypical representations of
  $C$. The $R(G)$-module structure on $R({\tH})(\pi)$ arises from its
  $R(H)$-module structure via the map $R(G) \to R(H)$, and since
  $R(H)$ is finite over $R(G)$ (\cite{SRRL}, Proposition 3.2), we are 
  reduced to showing that $R(\tH)(\pi)$ is finite as an $R(H)$-module.

  Now consider the restriction 
  \begin{equation}
    R(\tH) \to \prod_SR(S)
  \end{equation}
  where the $S$ runs through the conjugacy classes of Cartan subgroups
  of $\tH$ (conjugacy classes of such subgroups are finite in number
  and each one of the subgroups is closed, Abelian and contains the
  central subgroup 
  $C$). This map is injective (\cite{SRRL}, Proposition 1.2),
  whence $R(\tH)(\pi)$ is a subgroup of $\prod_SR(S)(\pi)$. Therefore
  it is enough to show that $R(S)(\pi)$ is finite as an $R(H)$-module
  for each $S$. The $R(H)$-module structure on $R(S)(\pi)$ arises from
  its structure of an $R(S/C)$-module via the map of representation
  rings induced by the inclusion $S/C \incl H$, and as $R(S/C)$ is
  finite over $R(H)$, it is enough to prove that that $R(S)(\pi)$ is
  finite over $R(S/C)$.

  We shall now show that $R(S)(\pi)$ is in fact a free $R(S/C)$-module
  with one generator. To prove this, recall that for a compact Abelian
  Lie group $A$, tensor product gives the set of irreducible
  representations $\hat{A}$ the structure of a finitely generated
  Abelian group, and that the representation ring of $A$ is given by
  the group ring $\Z[\hat{A}]$. Moreover, our exact sequence of
  compact Abelian groups
  \[
     1 \to C \to S \to S/C \to 1
  \]
  gives rise to an exact sequence
  \[
     1 \to \widehat{S/C} \to \hat{S} \to \hat{C} \to 1.
  \]
  From this it is clear that the summand $R(S)(\pi)$ of $R(S) =
  \Z[\hat{S}]$ is the subgroup freely generated by members of the
  coset of $\widehat{S/C}$ in $\hat{S}$ mapping to $\pi$ in $\hat{C}$,
  with the $R(S/C)$-module structure arising from the action of
  $\widehat{S/C}$ on the coset. Thus any representative of the coset 
  will form an $R(S/C)$-basis for $R(S)(\pi)$, and we are done.
\end{proof}
\end{proposition}

The following two lemmas point out further useful properties of the theories
$\proK_{\tG}^*(-)(\pi)\compl_{I_G}$.
\begin{lemma}
\label{rprokchg}
  Let $H$ be a closed subgroup of $G$, and let $X$ be a based $H$-CW complex.
Then there is a natural isomorphism of pro-$R(G)$-modules
\[
\rproK^*_\tG(G_+ \smashprod_H X)(\pi)\compl_{I_G} \isom
\rproK^*_\tH(X)(\pi)\compl_{I_H},
\]
where $\tH$ denotes the inverse image of $H$ in $\tG$.
\end{lemma}
\begin{proof}
Observe that the $H$-CW structure on $X$ gives rise to a $G$-CW
structure on $G_+ \smashprod_H X$, and that as $X_\alpha$ runs over
the finite $H$-CW subcomplexes of $X$, $G_+ \smashprod_H X_\alpha$
runs over the finite $G$-CW subcomplexes of $G_+ \smashprod_H X$. Now
the lemma follows from the $\tG$-equivariant isomorphism
\[
  \tG\, \smashprod_\tH X_\alpha \isom G \smashprod_H X_\alpha;
\]
from the change of groups isomorphism
 \[
  \rK^*_\tG(\tG \smashprod_\tH X_\alpha) \isom \rK^*_\tH(X_\alpha);
\]
from the compatibility of this isomorphism with the decomposition
\eqref{kdecomp}; and from the fact that the $I_G$-adic and $I_H$-adic
topologies on an $R(H)$-module coincide 
(see \cite{SRRL}, Corollary 3.9).
\end{proof}
\begin{lemma}
\label{freeprok}
Let $X$ be a free $G$-CW complex. Then there is a natural isomorphism
$\proK_\tG^*(X)(\pi)\compl_{I_G} \isom \proK_\tG^*(X)(\pi)$.
\end{lemma}
\begin{proof} (Compare with the proof of Proposition 4.3 in \cite{ASC}.)
Let $X_\alpha$ be a finite $G$-CW subcomplex of $X$, and let
$X_{\alpha,1},\ldots,X_{\alpha,k}$ be the $G$-CW 
subcomplexes of $X_\alpha$
such that $X_{\alpha,1}/G,\ldots,X_{\alpha,k}/G$ are the connected
components of $X_{\alpha}/G$. Since the action
of $G$ on $X$ is free, for each $i=1,\ldots,k$ we have an isomorphism
\[
   K_G(X_{\alpha,i}) \xrightarrow{\,\isom\,} K(X_{\alpha,i}/G).
\]
Pick a base point for $X_{\alpha,i} / G$. Then the diagram 
\[
\xymatrix{
  R(G) \ar[d] \ar[r] 
  & 
  K_G(X_{\alpha,i}) \ar[r]^{\isom\,\,} 
  & 
  K(X_{\alpha,i}/G) \ar[d] 
  \\
  \Z \ar@{=}[rr]
  & &
  \Z
}
\]
commutes, whence the composite of the maps in the top row sends $I_G$
into $\rK(X_{\alpha,i} / G)$. However, since $X_{\alpha,i} / G$ is a 
connected finite
CW-complex, the elements of $\rK(X_{\alpha,i} / G)$ are nilpotent. Because $R(G)$
is Noetherian, the ideal $I_G$ is finitely generated, and it follows that
for large enough $n$, the image of $I_G^n$ in $K_G(X_{\alpha,i})$
vanishes. As this happens for all $i=1,\ldots,k$,
the same is true of the image of $I_G^n$ in 
$K_G(X_{\alpha}) \isom \prod_{i=1}^k K_G(X_{\alpha,i})$.
Thus
\[
  I_G^n \cdot K_\tG^*(X_\alpha)(\pi) = 0
\]
for large $n$, and therefore
\begin{align*}
  \proK_\tG^*(X)(\pi)\compl_{I_G} 
      &= \{\quotb{K_\tG^*(X_\alpha)(\pi)}{I_G^n} \}_{\alpha,n}\\
      &\isom \{ K_\tG^*(X_\alpha)(\pi) \}_\alpha\\
      &= \proK_\tG^*(X)(\pi)
\end{align*}
as claimed.
\end{proof}

\begin{remark}
The main technical benefit of introducing the pro-group-valued
theories $\proK_\tG^*(-)(\pi)$ and
$\proK_\tG^*(-)(\pi)\compl_{I_G}$ is that they allow us to sidestep
problems with exactness that would otherwise complicate the proof 
of Theorem~\ref{mainthm}. The
source of these problems is the failure of inverse limits to preserve
exactness, as well as the failure of completion to be exact for
non-finitely generated modules. The idea of using pro-groups to prove
the completion theorem goes back to the original paper of Atiyah and Segal
\cite{ASC}.
\end{remark}

\section[Twistings arising from a central extension]{The case of a twisting arising from a graded central extension}
\label{sc}

In this section we will prove Theorem \ref{mainthm} in the case where
the twisting $\tau$ arises from a central extension in the way
explained in \cite{FHT1}. That is, we will prove the following.
\begin{theorem}
\label{babythm}
Let $X$ be a finite $G$-CW complex, where $G$ is a compact Lie group.
 Then the
projection $\pi : EG \times X \to X$ induces an isomorphism
\begin{displaymath}
  K^{\tau+*}_G(X)\compl_{I_G} \xrightarrow{\,\isom\,} K^{\pi^*(\tau)+*}_G(EG \times X)
\end{displaymath}
for any twisting $\tau$ arising from a graded central extension 
\[
    1\to \T \to G^{\tau} \to G \to 1,\qquad \epsilon : G \to \Z/2.
\]
\end{theorem}

Our argument for proving  Theorem \ref{babythm} is 
closely based on the one Adams, Haeberly, Jackowski and May present 
for proving a generalization of the Atiyah--Segal completion theorem 
in the untwisted case \cite{AHJM1}. Their argument in turn builds on 
ideas due to Carlsson \cite{Csbrc}.
As before, let
\[
   1 \to C \to \tG \to G \to 1 
\]
be a central extension of a compact Lie group $G$ by a compact 
commutative Lie
group $C$,
and let $\pi$ be an irreducible representation of $C$.
We shall derive Theorem \ref{babythm} from the following result.
\begin{theorem}
  \label{babythm2}
Suppose $X_1$ and $X_2$ are $G$-spaces, and let $f : X_1 \to X_2$ be a
$G$-equivariant map which is a non-equivariant homotopy equivalence.
Then the map
\[
 f^*: \proK^*_\tG(X_2)(\pi)\compl_{I_G} \to \proK^*_\tG(X_1)(\pi)\compl_{I_G}
\]
is an isomorphism.
\end{theorem}
Before proving Theorem \ref{babythm2}, we explain how
it implies Theorem \ref{babythm}.
\begin{proof}[Proof of Theorem \ref{babythm} assuming Theorem \ref{babythm2}]
Let $Z$ be a finite $G$-CW complex. By Theorem \ref{babythm2}, the
projection map $\pi : EG \times Z \to Z$ induces an isomorphism
\begin{equation}
\label{eqiso}
 \proK^*_{G^\tau}(Z)(1)\compl_{I_G} \xrightarrow[\isom]{\;\,\,\pi^*}
 \proK^*_{G^\tau}(EG \times Z)(1)\compl_{I_G}.
\end{equation}
Since $Z$ is finite, we have
\begin{equation}
\label{eqz}
  \begin{split}
  \proK_{G^\tau}^*(Z)(1)\compl_{I_G} &=
  \{\quotb{K_{G^\tau}^*(Z_\alpha)(1)}{I_G^n} \}_{\alpha,n}\\ &=
  \{\quotb{K_{G^\tau}^*(Z)(1)}{I_G^n} \}_n.
  \end{split}
\end{equation}
Fix a model for $EG$ which is a countable ascending union of
finite $G$-CW subcomplexes $EG_k$, $k \geq 1$; for example, we could
take $EG$ to be the iterated join construction of Milnor and take $EG_k$
to be the $k$-fold join of $G$ with itself. Then Lemma~\ref{freeprok}
and the finiteness of $Z$ imply that
\begin{equation}
\label{eqt}
  \begin{split}
      \proK_{G^\tau}^*(EG \times Z)(1)\compl_{I_G} &=
      \proK_{G^\tau}^*(EG \times Z)(1) \\ 
 &= \{ K_{G^\tau}^*(EG_k \times Z)(1) \}_k.
  \end{split}
\end{equation}
Thus applying the limit functor taking pro-$R(G)$-modules to
$R(G)$-modules to the isomorphism \eqref{eqiso} gives us an isomorphism 
\begin{equation}
\label{eqiso2}
K_{G^\tau}^*(Z)(1)\compl_{I_G} \xrightarrow[\isom]{\quad\,\,\pi^*\quad}
\invlim_{k} K_{G^\tau}^*(EG_k \times Z)(1).   
\end{equation}
Using \eqref{eqt}, \eqref{eqiso} and \eqref{eqz}, we see that inverse
system $\{K_{G^\tau}^*(EG_k\times Z)(1) \}_k$ is equivalent to one
that satisfies the Mittag--Leffler condition, whence the $\lim^1$
error terms vanish and the codomain in \eqref{eqiso2} is isomorphic to
$K_{G^{\tau}}^*(EG \times Z)(1)$.  Thus for any finite $G$-CW complex
$Z$, we have a natural isomorphism
\[
K_{G^\tau}^*(Z)(1)\compl_{I_G} \xrightarrow[\isom]{\quad\,\,\pi^*\quad}
K_{G^{\tau}}^*(EG \times Z)(1).
\]
Suppose now $Z$ is a based finite $G$-CW complex. Then from the diagram
\[
\xymatrix @C-0.5pc { 
0 \ar[r] & \rK_{G^\tau}^*(Z)(1)\compl_{I_G} \ar[r] \ar@{.>}[d] &
K_{G^\tau}^*(Z)(1)\compl_{I_G} \ar[r] \ar[d]^{\isom} &
K_{G^\tau}^*(\pt)(1)\compl_{I_G} \ar[r] \ar[d]^{\isom} & 0 \\ 
0 \ar[r] & \rK_{G^\tau}^*(EG_+ \smashprod Z)(1) \ar[r] &
K_{G^\tau}^*(EG \times Z)(1) \ar[r] & K_{G^\tau}^*(EG)(1) \ar[r] & 0
}
\]
we see that there is an induced isomorphism
\[
\rK_{G^\tau}^*(Z)(1)\compl_{I_G}
\xrightarrow[\isom]{\quad\,\,\pi^*\quad} \rK_{G^\tau}^*(EG_+ \smashprod Z)(1).
\]
The claim now follows by taking $Z$ to be the space $X_+ \smashprod
S^1(\epsilon)$ and applying Proposition~\ref{propceiso}.
\end{proof}
The rest of this section is dedicated to the proof of Theorem
\ref{babythm2}.  Let $\{V_i\}_{i \in I}$ be a set of representatives
for the isomorphism classes of the non-trivial irreducible complex
representations of $G$. Then $I$ is countable, the fixed-point 
subspace $V_i^G$ is zero for each
$i\in I$, and for every proper closed subgroup $H$ of $G$, there is some $i
\in I$ such that $V_i^H \neq 0$. Let $U$ be the direct sum of
countably infinite number of copies of $\bigoplus_{i\in I} V_i$, and let
\[ 
 Y = \colim_{V \subset U} S^V
\]
where the colimit is over all finite-dimensional $G$-subspaces of $U$ and
$S^V$ denotes the one-point compactification of $V$. 
Pick a $G$-invariant inner product on $U$, and observe that $Y^G$ is $S^0$.
\begin{lemma}
\label{lemmayhc}
The space $Y$ is $H$-equivariantly contractible for any proper closed
subgroup $H$ of $G$.
\end{lemma}
\begin{proof}
Since $Y$ has the structure of an $H$-CW complex, it is enough to show
that the fixed point set $Y^K$ is weakly equivalent to a point for any
subgroup $K$ of $H$. Given any finite-dimensional $G$-subspace $V
\subset U$, we can find a finite-dimensional $G$-subspace $W \subset
U$ such that $V \subset W$ and $(W-V)^K \neq 0$, where $W - V$ denotes
the orthogonal complement of $V$ in $W$.  But then the inclusion $S^V
\incl S^W$ is $K$-equivariantly null-homotopic, whence the map
$(S^V)^K \incl (S^W)^K$ is null-homotopic. Since $Y^K$ is given by the
union
\[
    Y^K = \colim_{V\subset U} (S^V)^K,
\]
the claim follows.
\end{proof}
\begin{lemma}
\label{lemmaypz}
  The pro-$R(G)$-module $\rproK^*_\tG(Y)(\pi)\compl_{I_G}$ is pro-zero.
\end{lemma}
\begin{proof}
For a finite-dimensional $G$-subspace $V\subset U$, let 
\[
 \lambda_V \in \rK_G(S^V) = \rK_\tG(S^V)(0) \subset \rK_\tG(S^V)
\]
denote the equivariant Bott class, where ``0'' refers to the trivial
representation of $C$. Then by Bott periodicity, each element of
$\rK^*_\tG(S^V)(\pi)$ is uniquely expressible as a product $x\lambda_V$, where
$x \in \rK_\tG^*(S^0)(\pi)$. Suppose $W \supset V$.
From the diagram
\[
\xymatrix{
 \rK_\tG^*(S^{W-V})(\pi) 
 \ar[r] 
 \ar[d]^{\smashprod \lambda_V}_{\isom} 
 & 
 \rK^*_\tG(S^0)(\pi) 
 \ar[d]^{\smashprod \lambda_V}_{\isom} 
 \\
 \rK_\tG^*(S^W)(\pi) \ar[r]  
 & 
 \rK^*_\tG(S^V)(\pi) 
}
\]
it follows that the map 
\[
 \rK_\tG^*(S^W)(\pi) \to \rK^*_\tG(S^V)(\pi) 
\]
sends the element $x\lambda_W = x\lambda_{W-V}\lambda_V$ to $x\chi^{\,}_{W-V}\lambda_V$, where
$\chi^{\,}_{W-V}$ denotes the image of $\lambda_{W-V}$ under the map
\[
   \rK_G(S^{W-V}) \to \rK_G(S^0)
\]
induced by the inclusion $S^0 \incl S^{W-V}$. Since this map is 
non-equivariantly null-homotopic, it follows from the diagram
\[
\xymatrix{
 \rK_G(S^{W-V}) \ar[d] \ar[r] & \rK_G(S^0) \ar[d] \ar@{=}[r] & R(G) \ar[d] \\
 \rK(S^{W-V})  \ar[r] & \rK(S^0)  \ar@{=}[r] & \Z
}
\]
that $\chi^{\,}_{W-V} \in I_G$. Thus if we choose $W \subset U$ so that it
is the direct sum of $V$ with $n$ $G$-invariant subspaces of $U$, then 
the map 
\[
 \quotb{\rK_G^*(S^W)(\pi)}{I_G^n} \to  \quotb{\rK_G^*(S^V)(\pi)}{I_G^n} 
\]
is zero. It follows that for any fixed $n$ the pro-$R(G)$-module 
\[
      \{ \quotb{\rK_\tG(S^V)}{I_G^n} \}_V
\]
is pro-zero, and therefore so is
\begin{equation*}
  \begin{split}
    \rproK_\tG(Y)(\pi)\compl_{I_G} &= \{ \quotb{\rK_\tG(S^V)}{I_G^n} \}_{n,V} \\
        &= \invlim_n \{ \quotb{\rK_\tG(S^V)}{I_G^n} \}_V
  \end{split}
\end{equation*}
where the inverse limit is taken in the category of pro-$R(G)$-modules.
\end{proof}
We are now ready to prove Theorem \ref{babythm2}.
\begin{proof}[Proof of Theorem  \ref{babythm2}]
It is enough to prove that $\rproK_\tG^*(Z)(\pi)\compl_{I_G}$
is pro-zero when $Z$ is a non-equivariantly contractible $G$-space;
the claim then follows by taking $Z$ to be the mapping cone of $f$. We
shall show that $\rproK_\tG^*(Z)(\pi)\compl_{I_G} = 0$ for such $Z$ by
induction on the subgroups of $G$, making use of the fact that any
strictly descending chain of closed subgroups of a Lie group is of
finite length.

To start the induction, we observe that in the case $G = \{e\}$ 
the claim follows from the assumption
that $Z$ is non-equivariantly contractible. Assume inductively that
\[
  \rproK_\tH^*(Z)(\pi)\compl_{I_H} = 0
\]
for all proper closed subgroups $H$ of $G$; here as before $\tH$ denotes 
the inverse image of $H$ in $\tG$. The inclusion of the fixed-point set
$Y^G = S^0$ into $Y$ gives a cofiber sequence 
\[
   S^0 \to Y \to Y/S^0
\]
whence we have a cofiber sequence
\[
   Z \to Z\smashprod Y \to Z \smashprod (Y/S^0).
\]
Thus to show that $\rproK_\tG^*(Z)(\pi)\compl_{I_G} = 0$, it is enough
to show that
\[ 
\rproK_\tG^*(Z\smashprod Y )(\pi)\compl_{I_G} = 0 
\] 
and 
\[
\rproK_\tG^*(Z \smashprod (Y/S^0))(\pi)\compl_{I_G} = 0.
\]

Let us first show that $\rproK_\tG^*(Z\smashprod Y )(\pi)\compl_{I_G} = 0$;
we claim that in fact \[\rproK_\tG^*(W\smashprod Y )(\pi)\compl_{I_G} = 0 \]
for any based $G$-CW complex $W$. Observing that 
\[
 \rproK_\tG^*(W \smashprod Y)(\pi)\compl_{I_G} = \invlim_\alpha
 \rproK_\tG^*(W_\alpha \smashprod Y)(\pi)\compl_{I_G}
\]
where $W_\alpha$ runs through all finite $G$-CW complexes of $W$, we
see that it is enough to consider the case where $W$ is
finite. Filtering $W$ by skeleta and working inductively reduces us to
the case where $W$ is of the form $G/H_+ \smashprod S^n$ for some $n$
and some closed subgroup $H$ of $G$, and using the suspension axiom
further reduces us to the case $W = G/H_+$. But now in the case $H =
G$ the claim follows from Lemma~\ref{lemmaypz}; and in the case
$H\lneq G$, it follows from the change of groups isomorphism (Lemma
\ref{rprokchg})
\[
   \rproK_\tG^*(G/H_+ \smashprod Y)(\pi)\compl_{I_G} 
  \isom  \rproK_\tH^*(Y)(\pi)\compl_{I_H}
\]
together with Lemma~\ref{lemmayhc}.

It remains to show that $\rproK_\tG^*(Z \smashprod
(Y/S^0))(\pi)\compl_{I_G} = 0$.  We shall show that in fact
\[\rproK_\tG^*(Z \smashprod W)(\pi)\compl_{I_G} = 0 \] for
any based $G$-CW complex $W$ such that $W^G$ is a point.  Arguing as
above, we see that it is enough to consider $W$ of the form $W
= G/H_+$, where $H$ now has to be a proper closed subgroup of $G$. But
in this case the claim follows from the change of groups isomorphism
(Lemma~\ref{rprokchg})
\[
 \rproK_\tG^*(Z \smashprod G/H_+)(\pi)\compl_{I_G} =
 \rproK_\tH^*(Z)(\pi)\compl_{I_H} 
\]
and the inductive assumption.
\end{proof}

\section{The general case}
\label{gc}
In this section we finally prove Theorem \ref{mainthm} in full
generality. We shall do so by considering successively more 
general spaces, starting with the case $X = \pt$ and proceeding by
change of groups and Mayer--Vietoris arguments. Since in general
completion is exact only for finitely generated modules, along the way
we check that the twisted $K$-groups that enter the Mayer--Vietoris
sequences are finitely generated over $R(G)$.

\begin{lemma}
\label{lemmapt}
Theorem \ref{mainthm} holds and $K_G^{\tau+*}(X)$ is
finitely generated over $R(G)$ when $X = \pt$.
\end{lemma}
\begin{proof}
  By \cite{FHT1}, Example 2.29, any twisting of a point arises from a 
graded central extension. Thus Theorem \ref{babythm} shows that 
Theorem \ref{mainthm} holds in this case. The claim about finite generation 
follows from Proposition \ref{propceiso} and the proof of Proposition
\ref{propexact}.
\end{proof} 

\begin{lemma}
\label{lemmaorbit}
Theorem \ref{mainthm} holds and $K_G^{\tau+*}(X)$ is finitely
generated over $R(G)$ when $X = G/H$, where $H$ is a closed subgroup
of $G$.
\end{lemma}
\begin{proof}
  Notice that $G/H = G \times_H \pt$ and that $EG \times G/H = G \times_H EG$.
For any $H$-space $Z$, we have a natural local equivalence of topological 
groupoids 
\[
  Z \modmod H \to G \times_H Z \modmod  G
\]
giving rise to a natural change of groups isomorphism
\begin{equation}
\label{twchg}
   K^{\tau + *}_G(G \times_H Z) \xrightarrow{\quad\isom\quad} K^{\tau + *}_H(Z).  
\end{equation}
Consider the diagram
\[
\xymatrix{ 
K^{\tau+*}_G(G/H)\compl_{I_G} \ar[d] \ar[r] &
  K^{\tau+*}_H(\pt)\compl_{I_H} \ar[d] \\
K^{\pi^*(\tau)+*}_G(EG \times G/H) \ar[r] & K_H^{\pi^*(\tau)+*}(EG)
}
\]
Here the bottom row is a change of groups isomorphism as in \eqref{twchg};
the top row is an isomorphism because of the isomorphism \eqref{twchg} and
the fact that $I_H$-adic and $I_G$-adic completions of an $R(H)$-module agree
(see \cite{SRRL}, Corollary 3.9); and the vertical map on the right is an
isomorphism by Lemma~\ref{lemmapt} and the observation that $EG$ is a model
for $EH$. Thus the map on the left is also an isomorphism, which shows that 
Theorem \ref{mainthm} holds in this case. To see that $K^{\tau+*}_G(G/H)$
is finitely generated as an $R(G)$-module, observe that the isomorphism 
 \eqref{twchg} and Lemma~\ref{lemmapt} imply that it is finitely 
generated over $R(H)$. The claim now follows from the fact that
$R(H)$ is finite over $R(G)$ (see \cite{SRRL}, Proposition 3.2).
\end{proof}
\begin{lemma}
\label{lemmasph}
Theorem \ref{mainthm} holds and $K_G^{\tau+*}(X)$ is finitely
generated over $R(G)$ when $X$ is of the form $X = G/H \times S^n$, $n\geq0$.
\end{lemma}
\begin{proof}
The case where $n=0$ follows from Lemma~\ref{lemmaorbit} and the axiom
of disjoint unions. For $n > 0$, the claim follows inductively from
the Mayer--Vietoris sequences arising from the
decomposition of $S^n$ into upper and lower hemispheres $S^n_+$ and
$S^n_-$. Lemma~\ref{lemmaorbit} and the inductive 
assumption imply that all groups in the Mayer--Vietoris sequence
  \begin{multline}
\label{mvseq}
  \cdots \to K^{\tau+*}_G(G/H \times S^n) \to \\ \to K^{\tau+*}_G(G/H \times
  S^n_+) \oplus K^{\tau+*}_G(G/H \times S^n_+) \to \\ \to
  K^{\tau+*}_G(G/H \times (S^n_+ \cap S^n_-)) \to \cdots
  \end{multline}
except $K^{\tau+*}_G(G/H \times S^n)$ are finitely generated
over $R(G)$, whence the remaining group 
$K^{\tau+*}_G(G/H \times S^n)$ must also be finitely
generated, as claimed. It follows that the sequence obtained from
\eqref{mvseq} by completion with respect to the augmentation ideal
$I_G$ is exact.  Now the claim that Theorem \ref{mainthm} holds for
the space $G/H \times S^n$ follows from Lemma~\ref{lemmaorbit} and the
inductive assumption by comparing the completed sequence to the
Mayer--Vietoris sequence of the pair 
\[
(EG \times G/H \times S^n_+,EG \times G/H \times S^n_-)
\]
and applying the five lemma.
\end{proof}

Theorem \ref{mainthm} is now contained in the following theorem.

\begin{theorem}
Theorem \ref{mainthm} holds and $K_G^{\tau+*}(X)$ is finitely
generated over $R(G)$ for any finite $G$-CW complex $X$.
\end{theorem}

\begin{proof}
We proceed by induction on the number of cells in $X$. If $X$ has
no cells, that is, if $X$ is the empty $G$-space, the claim holds
trivially.  Assume inductively that the claim holds for the
space $X$, and consider the space $Y = X \cup_f (G/H \times
D^n)$, where $f : G/H \times S^{n-1} \to X$ is an attaching
map. Denote
\[
 D^n(r) = \{ x\in \R^n : |x| \leq r \},
\]
and let \[Y_1 = X \cup_f (G/H \times (D^n - D^n(1/3))) \subset Y\] and
\[Y_2 = D^n(2/3) \subset Y. \] 
By Lemma~\ref{lemmaorbit}, Lemma~\ref{lemmasph} and the inductive
assumption, in the Mayer--Vietoris sequence
\begin{equation*}
\xymatrix@1@C=1.2em{
  \cdots \ar[r] &  
  K_G^{\tau+*}(Y) \ar[r] &
  K_G^{\tau+*}(Y_1) \oplus  K_G^{\tau+*}(Y_2)\ar[r] &
  K_G^{\tau+*}(Y_1 \cap Y_2) \ar[r] &
  \cdots
}
\end{equation*}
all groups except possibly $K_G^{\tau+*}(Y)$ are finitely generated
over $R(G)$.  It follows that $K_G^{\tau+*}(Y)$ is also finitely
generated, as claimed. We conclude that the top row in the following diagram
of Mayer--Vietoris sequences is exact.
\[
\def\objectstyle{\scriptstyle}\def\labelstyle{\scriptstyle}
\xymatrix @C=1.2em { \cdots  \ar[r] & K^{\tau+*}_G(Y)\compl_{I_G} \ar[r] \ar[d] &
  K^{\tau+*}_G(Y_1)\compl_{I_G} \oplus K^{\tau+*}_G(Y_2)\compl_{I_G}
  \ar[r] \ar[d]^{\isom} & K^{\tau+*}_G(Y_1 \cap Y_2)\compl_{I_G} \ar[d]^{\isom}
  \ar[r] & \cdots \\ \cdots \ar[r] & K^{\tau+*}_G(EG\times Y) \ar[r] &
  K^{\tau+*}_G(EG\times Y_1) \oplus K^{\tau+*}_G(EG\times Y_2) \ar[r]
  & K^{\tau+*}_G(EG \times (Y_1 \cap Y_2)) \ar[r] & \cdots }
\]
In the diagram, the vertical map on the right is an isomorphism by
Lemma~\ref{lemmasph} and the map in the middle is an isomorphism by
Lemma~\ref{lemmaorbit} and the inductive assumption. Thus the map on
the left is an isomorphism by the five lemma, showing that
Theorem~\ref{mainthm} holds for the space $Y$.
\end{proof}

\bibliographystyle{plain}

\bibliography{completion-theorem}

\end{document}